\documentclass[12pt,psamsfonts]{article}
\usepackage{amsmath}
\usepackage{amssymb}
\usepackage{amsthm}

\newtheorem{defn0}{Definition}[section]
\newtheorem{prop0}[defn0]{Proposition}
\newtheorem{thm0}[defn0]{Theorem}
\newtheorem{lemma0}[defn0]{Lemma}
\newtheorem{corollary0}[defn0]{Corollary}
\newtheorem{example0}[defn0]{Example}
\newtheorem{remark0}[defn0]{Remark}
\newtheorem{conjecture0}[defn0]{Conjecture}

\newenvironment{proposition}{\bigskip \begin{prop0}}{\end{prop0}}
\newenvironment{theorem}{\bigskip \begin{thm0}}{\end{thm0}}
\newenvironment{lemma}{\bigskip \begin{lemma0}}{\end{lemma0}}
\newenvironment{corollary}{\bigskip \begin{corollary0}}{\end{corollary0}}

\newenvironment{remark}{ \begin{remark0}\rm}{\end{remark0}}

\newcommand{\propref}[1]{Proposition~\ref{#1}}
\newcommand{\thmref}[1]{Theorem~\ref{#1}}
\newcommand{\lemref}[1]{Lemma~\ref{#1}}
\newcommand{\corref}[1]{Corollary~\ref{#1}}

\newcommand{\secref}[1]{Section~\ref{#1}}

\def\cocoa
{\mbox{\rm C\kern-.13em o\kern-.07
em C\kern-.13em o\kern-.15em A}}

\def\max{{\mathfrak{m}}}                   
\def\res{{\bf k}}                   

\newcommand{\m}{\mathfrak{m}}
\newcommand{\n}{\mathfrak{n}}

\newcommand{\la}{\lambda}

\title{  \bf \huge Structure theorems for certain  Gorenstein  ideals.
}
\author{\large   Juan Elias
\thanks{Partially supported by  MTM2007-67493}
\and \large Giuseppe Valla \
\thanks{Partially supported by the Consiglio
Nazionale delle Ricerche (CNR)}}

\date{June 6, 2007}

\begin{document}

\maketitle


\bigskip
\section{Introduction.}

Let $I$ be an ideal in the regular local ring $(R,\n)$ such that $I\subseteq \n^2$  and let
$$A:=R/I, \ \ \m:=\n/I,\ \ \res:=R/\n=A/\m.$$ Let  $d=\dim(A)$ be the dimension, $e$ the multiplicity
 and $h=v(\m)-d$  the embedding codimension of $A.$
 We  assume that $\res$ is a characteristic zero field (see the comment after
 \propref{lean}).

 A classical problem in the theory  of local rings is the determination of the minimal number of
 generators $v(I):=\dim_k(I/\n I)$ of the ideal $I$ under certain restrictions on the numerical
 characters of $A$. For example, by a classical theorem of Abhyankar, we know that $e\ge h+1,$
 and if the equality $e=h+1$ holds we say that $A$ has minimal multiplicity and we know that
 $ v(I)={{h+1}\choose{2}}.$

 In a sequence of papers Rosales and Garc\'{\i}a-S\'{a}nchez
  proved the following   results in the
 case $A$ is the one dimensional local domain  corresponding to a monomial curve in the affine space,
 see,  \cite{Ros96b}, \cite{Ros96a}, \cite{RGS98}.
 By very hard computations related to the numerical semigroup of the curve, they were able to prove that

\noindent
If $h+2\le e \le h+3$, then
\begin{equation}\label {R1}
{{h+2}\choose{2}}-e\le v(I)\le {{h+1}\choose{2}}.
\end{equation}
\noindent
If $h+2\le e \le h+4 $ and $A$ is Gorenstein, then
\begin{equation}\label {R2}
v(I)={{h+1}\choose{2}}-1.
\end{equation}
We remark that the monomial curve $\{t^8:t^{10}:t^{12}:t^{15}\}$
shows that (\ref{R2}) does not hold if $e=h+5$, see \cite{RGS98}.

On the other hand, the  monomial curve $\{t^7:t^8:t^{10}:t^{19}\}$
shows that the upper bound in (\ref{R1}) does not hold if $e=h+4$.
In the same paper it is asked whether it is true that, with $e=h+4$,
one has
\begin{equation}\label {R3}
{{h+2}\choose{2}}-e={{h+1}\choose{2}}-3\le v(I)\le
{{h+1}\choose{2}}+1.
\end{equation}
\noindent

A first motivation for our paper was to understand  these results and
to extend them to the general case of a local Cohen-Macaulay ring of any dimension.

\bigskip
A sharp upper bound for the minimal number of generators of a
perfect ideal $I$ in a regular local ring $R,$ has been given in
\cite{ERV91} in terms of the multiplicity $e$  and of the
codimension $h$ of $R/I.$ The bound is $$ v(I)\le
\binom{h+t-1}{t}-r+r^{<t>},$$ where the meaning of $r,t$ and
$r^{<t>}$ will be explained in the \secref{preli}. In the same
section we will also prove that $$ {{h+2}\choose {2}}-e\le v(I)$$
holds for every perfect codimension $h$ ideal $I$ in a regular local
ring $R,$ see Proposition (\ref{prop}). Further we will see  how
these bounds extend (\ref{R1}) to a considerable extent and
positively answer question (\ref{R3}) in a very general setting.

\vskip 2mm As for (\ref{R2}), the problem is much harder. We have a
Gorenstein local ring $(A=R/I,\m=\n/I)$ of codimension $h$ and
multiplicity $h+2\le e \le h+4$ and we want to determine the minimal
number of generators of $I$. It is easy to see that we may  assume
that $A=R/I$ is artinian; since $A$ is Gorenstein, the possible
Hilbert function of $R/I$ are
$$(1,h,1),(1,h,1,1),(1,h,2,1),(1,h,1,1,1),$$  so that, in any case,
$v(\m^2)\le 2.$

Following Sally (see \cite{Sal79c}), we say that an Artinian local
ring $(A,\m),$ not necessarily Gorenstein, is {\bf stretched} if
$v(\m^2)=1.$ We call {\bf Almost stretched} an Artinian local ring
such that $v(\m^2)=2.$

With this notation, we strongly extend (\ref{R2}) if we can prove
that if $R/I$ is Gorenstein, stretched or almost stretched of
multiplicity $e$ and codimension $h,$  then $v(I)=\binom{h+1}{2}-1.$

By the classical theorem of Macaulay on the shape of the Hilbert
Function of a standard graded algebra,   the Hilbert function of $A$
is given by:

\begin{center} \begin{tabular}{|c|c|c|c|c|c|}
0 & 1  & 2 & \dots  & s  & s+1  \\ \hline
1  & h  & 1 & \dots & 1  &  0  \\
\end{tabular} \end{center}
with $(s\ge 2)$ if $A$ is stretched, or by

\begin{center} \begin{tabular}{|c|c|c|c|c|c|c|c|c|c|c|c|}
0 & 1  & 2 & \dots & t & t+1  &\dots & s  & s+1  \\ \hline
1  & h  & 2 & \dots & 2 & 1 &  \dots & 1  &  0  \\
\end{tabular} \end{center}
with $s\ge t \ge 2,$ if $A$ is almost stretched

\vskip 2mm

The particular shape of the Hilbert function can be used to prove
that

$$
  \binom{h+1}{2}-1 \le v(I)\le \binom{h+1}{2} \qquad \text{if $A$ is stretched,\phantom{ almost}}
  $$

  $$  \binom{h+1}{2}-2 \le v(I)\le \binom{h+1}{2} \qquad \text{if $A$ is almost stretched}.
$$

The case of stretched Artinian Gorenstein local ring has been
studied by J. Sally in \cite{Sal79c}  where she was able to prove a
structure theorem for the corresponding ideals, see also
\cite{RosV01}. We extend this result    to the case of  stretched
Artinian local rings of any Cohen-Macaulay type. But an   unexpected
and  deeper result which we will prove in this paper, is a structure
theorem for {\it any almost stretched Gorenstein local rings}.

\vskip 2mm These  results  are proved in Section 3 and 4, Theorem
\ref{first} and  Theorem \ref{main}, respectively.

\vskip 2mm Of course, as a consequence, we get even more of what we
wanted, namely:

\bigskip
 $ v(I)=\binom{h+1}{2}-1$  if $A$ is stretched and  $\tau(A)<h,$ while $ v(I)=\binom{h+1}{2}$
 otherwise;

\bigskip $ v(I)=\binom{h+1}{2}-1$ if $A$ is almost stretched and
Gorenstein.

\vskip 10mm Another motivation for our paper came from a recent work
by Casnati and Notari (see \cite{CN07}).
Let $\mathcal
Hilb_{p(t)}(\mathbb P ^n_k)$ denote the Hilbert scheme parametrizing
closed subschemes in $\mathbb P ^n_k$ with given Hilbert polynomial
$p(t)\in \mathbb Q[t].$

The case $\deg (p(t))=0$ is often  problematic. Since it is known
that  any zero-dimensional Gorenstein  scheme of degree d can be
embedded as an arithmetically Gorenstein non-degenerate subscheme in
$\mathbb P ^{d-2}_k,$
 it is natural to study the open locus
 $$\mathcal Hilb_d^{aG}(\mathbb P ^{d-2}_k) \subseteq \mathcal Hilb_d(\mathbb P ^{d-2}_k).$$

 The scheme $\mathcal Hilb_d^{aG}(\mathbb P ^{d-2}_k) $ has a natural stratification
 which reduce the problem to understand the intrinsic structure of Artinian Gorenstein
 $\res$-algebras of degree $d$. Since such an algebra is the direct sum of local,
 Artinian, Gorenstein $\res$-algebras of degree at most $d$, it is natural to begin
 with the inspection of these elementary bricks.

 If $d=6,$ the bricks are all given by stretched local rings, save for the case of
 Hilbert function $(1,2,2,1)$ which is almost stretched and was studied deeply by Casnati and Notari.

  If we want to extend the above results to the case $d\ge 7$, the first step is to
  study the intrinsic structure of Artinian Gorenstein local algebras  with multiplicity $7$.
Since the  Hilbert  function $(1,2,3,1)$ is not allowed,
 an Artinian Gorenstein ring $(A,\m)$ with multiplicity $7$ is stretched or almost stretched.
See \cite{Iar87b} for more results on the classification of Artin
algebras.

Hence, the structure theorems we will prove in the next sections
will  give light to these questions too.

It is clear that the best would be to have a classification up to
isomorphisms of artinian Gorenstein k-algebras of a given Hilbert
function, at least in the almost stretched case. We approach this
very difficult problem in the last part of the paper, where we give
a classification of Artinian complete intersection local k-algebras
with Hilbert function $(1,2,2,2,1,1,1).$  This example is
significant because the parameter space has a one-dimensional
component.


\section{Upper and lower bounds for $v(I).$}

\label{preli}

Let $(R,\n)$ be a regular local ring, $I$ an ideal in $R$. Let us
assume that $(A=R/I,\m=\n/I)$ has dimension $d$, embedding
codimension $h$ and multiplicity $e.$
We denote by $H_A$ the Hilbert function of $A$
$$
H_A(n):=dim_{\res}\left(\frac{\max^n}{\max^{n+1}}\right)
$$
$n\ge 0$. The socle degree of an Artin ring $A$ is the last integer
$s=s(A)$ such that $H_A(s)\neq 0$; the Cohen-Macaulay type of $A$
is
$$
\tau(A):=dim_{\res}(0:\max).
$$

A sharp upper bound for $v(I)$ can be given by using the notion of
lex-segment ideal as in  \cite{ERV91}. We recall that the associated
graded ring of $A$ can be presented as  $gr_{\m}(A)=gr_{\n}(R)/I^*,$
where $I^*$ is the ideal generated by the $\n$-initial forms of $I$
in the polynomial ring $S=gr_{\n}(R).$ This implies that the Hilbert
Function of $A=R/I$ is the same as the Hilbert Function of the
standard graded algebra $S/I^*.$

A set of elements in $I$ whose $\n$-initial forms generate $I^*,$ is
called a {\it standard basis} of $I.$
  Since it is easy to see that   a standard basis is a basis,  we have the inequality $v(I)\le v(I^*).$

On the other hand, by a classical result of Macaulay, any
homogeneous ideal $P$ in the polynomial ring $S=k[X_1,\dots,X_n]$
has the following property: the number of minimal generators of $P$
is less than or equal to the  number of minimal generators of the
unique lex-segment  ideal $P_{lex},$ which has the same Hilbert
Function of $P.$

Hence, given the  ideal $I$ in the regular local ring $(R,\n)$ and
the corresponding lex-segment ideal $I_{lex}:=(I^*)_{lex}$ in
$S:=gr_{\n}(R),$  we have \begin{equation}\label{lex} v(I)\le
v(I^*)\le v(I_{lex}).\end{equation}

 More difficult is to get a bound only involving the multiplicity and the codimension.
 Namely one has to compare the number of generators of all the lex-segment ideals having
 the given multiplicity and codimension. This has been done  in \cite{ERV91}
  where the following bound has been proved.

 We need some more notations. If $n$ and $i$ are positive integers then $n$ can be uniquely
 written as $$n=\binom {n(i)} {i}+\binom {n(i-1)} {i-1}+\dots +\binom {n(j)} {j}$$
 where $n(i)>n(i-1)>\dots >n(j)\ge j \ge 1.$ This is called the $i$-binomial expansion of $n.$ We let
 $$n^{<i>}:=\binom {n(i)+1} {i+1}+\binom {n(i-1)+1} {i}+\dots +\binom {n(j)+1} {j+1}.$$

 Given two positive integers $e,h$ with $e\ge h+1$ we define $t$ as the unique integer such that
 $$\binom{h+t-1}{t-1}\le e < \binom{h+t}{t}$$ and $$r:=e-\binom{h+t-1}{t-1}.$$
 The main result in \cite{ERV91} shows that, for every perfect codimension $h$ ideal $I$
 in the regular local ring $R$ with $I\subseteq \n^2$ and $e(R/I)=e,$ we have
 \begin{equation}\label{ERV} v(I)\le \binom{h+t-1}{t}-r+r^{<t>}.\end{equation}
 For example if $h\ge 3$ and $e=h+2,$ then $t=2,$ $r=1$ and we get $v(I)\le \binom{h+1}{2}.$
 The same bound holds also for $e=h+3,$ see (\ref{R1}).

 Instead, if $e=h+4$ we get $t=2,$ $r=3$ and
 $$v(I)\le \binom{h+1}{2}-3+3^{<2>}=\binom{h+1}{2}-3+4=\binom{h+1}{2}+1,$$
  see (\ref{R3}). The same bound holds also for $e=h+5.$

  \vskip 2mm A lower bound for $v(I)$ follows from the following easy lemma.
  \begin{lemma} \label{art} Let $A=R/I$ be a local Artinian ring with multiplicity
  $e$ and embedding codimension $h.$ We assume that $I\subseteq \n^2.$
  Then we have $$ \binom{h+2}{2}-e\le  \binom{h+1}{2}-v(\m^2)\le v(I).$$
 \end{lemma}
 \begin{proof}
 It is clear that the Kernel of the epimorphism
 $$\n^2/\n^3\to \m^2/\m^3=(\n^2+I)/(\n^3+I)\to 0$$ is
 $(\n^3+I)/\n^3\cong I/(\n^3\cap I).$
 Since $I\n\subseteq \n^3\cap I,$ we get $$v(\n^2)-v(\m^2)={{h+1}\choose{2}}-v(\m^2)\le v(I).$$

 Notice that we  have $e=\sum_{i=0}^sv(\m^i)$, where $s$ is the socle degree of $A$,  so that
 $e\ge 1+h+v(\m^2)$ and $$\binom{h+2}{2}-e\le \binom{h+2}{2}-(1+h+v(\m^2))=\binom{h+1}{2}-v(\m^2).$$
\end{proof}

As a consequence of this lemma we get a lower bound for the number of
generators of perfect ideals in a regular local ring which, at least
for low multiplicity, seems to  be useful.

\begin{proposition}\label{prop}  Let Let $A=R/I$ be a local Cohen-Macaulay ring with
dimension $d$, multiplicity $e$ and embedding codimension $h.$ We assume that $I\subseteq \n^2.$
Then we have $$\binom{h+2}{2}-e\le v(I)\le \binom{h+t-1}{t}-r+r^{<t>}.$$
\end{proposition}
\begin{proof} Let $J=(x_1,...,x_d)$  be a maximal $\n$-superficial sequence for $A.$
Since $A$ is Cohen-Macaulay,  $x_1,...,x_d$ is a regular sequence modulo $I$ so that  $I\cap J=IJ.$
 Let $$\overline{I}=(I+J)/J,\ \ \ \overline{R}=R/J, \ \ \ \overline{A}=A/(x_1,...,x_d)A= \overline{R}/\overline{I},
 \ \ \ \overline{\m}=\m/J.$$
 Then we have $$v(\overline{I})=\dim_k(I+J/\n I+J)=\dim_k(I/\n I+I\cap J)=\dim_k(I/\n I)=v(I).$$
We know also that the multiplicity of $A$ is the same as the
multiplicity of the Artinian local ring $A/(x_1,...,x_d)A$. Finally  $I$ and
$\overline{I}$ share the same  embedding codimension because
$h=v(\m)-d=v(\overline{\m}).$  The lower bound now follows from
Lemma \ref{art}, while the upper bound is given by (\ref{ERV}).
\end{proof}

\vskip 3mm In the next section we are going to establish structure
theorems for stretched local rings and for almost stretched
Gorenstein local rings. One of the main ingredient will be the
following result which will be used several times later and  is
reminiscent of the {\bf lean basis} notion introduced by J. Sally in
\cite{Sal79c}.

In the proof of the following Proposition we need to know that if
the characteristic of $\res$ is $0$, then {\bf a Borel fixed monomial
ideal $K$ is strongly stable}. This means that $K$  satisfies the
following requirement: for any term $M\in K$ and any indeterminate
$X_j$ dividing $M$, we have $X_i(M/X_j )\in K$ for all $1\le i <j.$

\begin{proposition}\label{lean} Let $(A,\m)$ be an Artinian local ring of embedding dimension $h$ and socle degree $s$
such that the characteristic of the residue field $\res$ is $0$ and $v(\m^2)\le 2.$
Then we can find a minimal basis $x_1,\dots,x_h$ of $\m$ such that
$$\m^j=(x_h^j),\ \ \ \ j=2,\dots,s $$  if $A$ is stretched, while
$$\m^j=\begin{cases}
    (x_h^j,x_h^{j-1}x_{h-1}) & j=2,\dots,t \\ (x_h^j)  & j=t+1,\dots,s
\end{cases}$$ if $A$ is almost stretched.
\end{proposition}
\begin{proof} We prove the proposition in the case $A$ is almost stretched,
 because the other case is easier.
Let $\m=(a_1,\dots,a_h);$ we know that the Hilbert function of $A$
is the same as the Hilbert function of
$gr_{\m}(A)=k[\xi_1,\dots,\xi_h]=S/J$ where
$\xi_i:=\overline{a_i}\in \m/\m^2,$ $S=k[X_1,\dots,X_h]$,  and $J$
is an homogeneous ideal of $S.$ Further, the generic initial ideal
$gin(J)$ of $J$ is a Borel fixed monomial ideal which is then
strongly stable.

We claim that, after a suitable changing of coordinates in $S$,
which corresponds to a changing  of generators for the maximal ideal
$\m$ of $A,$ we may assume that  a basis for $S_j$ modulo $gin(J)_j$
is given by $X_h^j,X_h^{j-1}X_{h-1}$ for $ j=2,\dots,t,$ and by
$X_h^j$ for $ j=t+1,\dots,s .$

In order to prove the claim, we need only to remark that if a
monomial  ideal $K$ is strongly stable  and $K_j\not=S_j,$ then
$X_h^j\notin K_j$, and if $\dim_k(S_j/K_j)\ge 2,$ then also
$X_h^{j-1}X_{h-1} \notin K_j$.

Since $gin(J)$ is an initial ideal, the same monomials form a basis
also   for $S$ modulo $J.$  The conclusion follows because we have
for every $j\ge 0$ $$S_j/(J)_j=(\m^j/\m^{j+1}).$$

\end{proof}

Because of this Proposition, we will always assume in the paper that
{\bf the residue field $\res$ has characteristic zero}.

\begin{remark}
Notice  that if  the codimension is bigger than two, the argument used in the proof of
\propref{lean} is not true anymore. Take for example the ideals
$(X_1^2,X_1X_2,X_1X_3)$ and $(X_1^2,X_1X_2,X_2^2)$ which are
strongly stable of codimension three in
$k[X_1,X_2,X_3].$
\end{remark}

\bigskip
\section{Stretched local rings}
\label{SLC}
We recall that in \cite{Sal79c} J. Sally studied several properties of stretched
local rings and proved a structure theorem for stretched Artinian local rings
in the Gorenstein case. Here we extend the result to any Cohen-Macaulay type.

\begin{theorem}\label{first}   Let $I$ be an ideal in the regular local ring
$(R,\n)$ such that $I\subseteq \n^2$ and $A:=R/I$ is Artinian. Let $\m:=\n/I,$  $h:=v(\m)$
and $\tau$ the Cohen-Macaulay type of $A.$

\vskip 2mm
\noindent
$(1)$
If $A$ is stretched of socle degree $s$ and $\tau<h,$ then we can find a  basis
$\{x_1,\dots,x_h\}$ of $\n$ such that  $I$ is minimally generated by the elements
$\{x_ix_j\}_{1\le i<j\le h}$, $\{x_j^2\}_{2\le j \le \tau},$
$\{x_i^2-u_ix_1^s\}_{ \tau +1\le i \le h},$ where the $u_i$ are units in $R.$

\vskip 2mm \noindent
$(2)$ If $A$ is stretched of socle degree $s$ and $\tau=h,$ then we can find
 a  basis $\{x_1,x_2, \dots,x_h \}$ of $\n$ such that  $I$ is minimally
 generated by the elements $\{x_1x_j\}_{2\le j\le h}$, $\{x_ix_j\}_{2\le i  \le j\le h}$ and $x_1^{s+1}.$
\end{theorem}
\begin{proof}
By  Proposition \ref{lean}, we can find an element $y_1\in \m, y_1\notin \m^2$
such that $y_1^s\not=0$ and  $\m^j=(y_1^j)$ for  $ 2\le j \le s.$
We remark that this implies $y_1^j\notin \m^{j+1}$ for every $1\le j \le s.$

\begin{lemma} We have $$(0:\m)\bigcap \m^2=\m^s.$$
\end{lemma}
\begin{proof}
If $s=2$ there is nothing to prove, hence let $s\ge 3.$ If $a\in 0:\m$
and $a\in \m^2,$ then $a=y_1^2u$ and we get $0=y_1a=y_1^3u.$ Since $s\ge 3,$
this implies   $u\in \m,$ otherwise $y_1^3=0.$ Hence  $a\in \m^3;$
 going on in this way we get $a\in \m^s$ as wonted.
\end{proof}

Since  $y_1^s\in 0:\m$ and $y_1^s\not=0,$   we can find elements
$y_2,\dots,y_{\tau}\in \m$ such that  $\{y_1^s,y_2,\dots,y_{\tau}\}$ is a basis of the $\res$-vector space $0:\m.$

\begin{lemma}
The elements $y_1,y_2,\dots,y_{\tau}$ are part of a minimal basis of $\m.$
\end{lemma}
\begin{proof}
If $\sum_{i=1}^{\tau} \la_iy_i \in \m ^2,$  then $\la_1 \in \m$,
otherwise $y_1\in 0:\m+\m^2$ and $y_1^2\in \m^3,$ a contradiction.
Thus we get $$\sum_{i=2}^{\tau} \la_iy_i\in (0:\m)\bigcap \m^2=\m^s$$ and,
for some $t\in R,$  $\sum_{i=2}^{\tau} \la_iy_i+ty_1^s=0.$
 This  implies $\la_i \in \m$ for every $i,$ because  $\{y_2,\dots,y_{\tau},y_1^s\}$
 is a basis of the $\res=A/\m$ vector space $0:\m.$
\end{proof}

Of course  we can complete the set  $\{y_1,y_2,\dots,y_{\tau}\}$ to
a minimal basis of $\m,$ say  $\m=(y_1,y_2,\dots,y_{\tau},z_{\tau+1},\dots,z_h).$
Now, if $j \ge \tau +1,$ we have $y_1z_j\in \m^2,$ hence $y_1z_j=y_1^2t$
and $z_j-y_1t \in 0:y_1.$ By replacing $z_j$ with $z_j-y_1t$ in the minimal
 generators of $\m$, we may assume that $$\m=(y_1,y_2,\dots,y_{\tau},y_{\tau +1},\dots,y_h)$$
 with
 \begin{equation}
 \label{m}y_2,\dots,y_{\tau}\in 0:\m,\ \ \ \ \  y_{\tau +1},\dots,y_h\in 0:y_1.
 \end{equation}

\medskip
\noindent
{\bf Let us first consider the case $\tau<h.$}

\vskip 2mm \noindent If we choose $i$ and $j$ so that $\tau +1\le i\le j\le h,$
we have $$y_iy_j\m\subseteq y_i\m^2=y_i(y_1^2)=0.$$ Hence $y_iy_j\in (0:\m)\cap \m^2=\m^s=(y_1^s),$
and we can write $y_iy_j=u_{ij}y_1^s$ where   $u_{ij}\in \m$ if and only if $y_iy_j=0.$

If we let $J:=(y_{\tau +1},\dots,y_h),$ we may define an inner product in the
 $\res$-vector space $V:=J/J\m$ by letting $$<\overline{y_i},\overline{y_j}>:=\overline{u_{ij}}\in A/\m=\res.$$
 This is well defined.  Namely,  let $y_i=p_i+z_i$ with $p_i\in J$ and $z_i\in J\m;$
 since $J\subseteq 0:y_1,$ we get  $$y_iy_j-p_ip_j=(p_i+z_i)(p_j+z_j)-p_ip_j\in J\m^2=y_1^2J=0.$$
Since the characteristic of $\res$ is not two, the inner product can be diagonalized.
This means that the generators of $\m$ can be chosen to satisfy
 \begin{equation}
 \label{p} y_iy_j=0
 \end{equation}
  for every  $\tau +1\le i< j\le h.$
   This implies that   for every $\tau +1\le i\le h,$ we must have $y_i^2\not=0,$
   because, if $y_i^2=0,$ we would get $y_i\in 0:\m,$  a contradiction. Hence,
   for every $\tau +1\le i\le h$,  we will have
   \begin{equation}\label{b}
   y_i^2=u_iy_1^s
   \end{equation} with $u_i\notin \m.$

As a consequence  we can prove the first part of the   theorem.
Let $x_i\in \n$ such that $\overline{x_i}=y_i.$
From (\ref{m}), (\ref{p}) and (\ref{b}),  it is clear that all the  elements
$$
\{x_ix_j\}_{1\le i<j\le h},\ \
\{x_j^2\}_{2\le j \le \tau},\ \
\{x_i^2-u_ix_1^s\}_{ \tau +1\le i \le h},
$$ are in $I.$
Let $J$ be the ideal they generate; then $J\subseteq I$ so that
$H_{R/I}(n)\le H_{R/J}(n)$ for every $n\ge 0.$
We claim that we have equality above  for every $n\ge 0.$
Namely we have $$x_1^{s+1}=(u_h)^{-1}x_1x_h^2\in J$$ so that $I^*\supseteq J^*\supseteq K$
where $K$ is the ideal in $S= \res[X_1,\dots,X_h]$ generated by
$X_1^{s+1}$ and all  degree two monomials  except $X_1^2.$
Since the Hilbert function of $S/K$ is the same as the Hilbert function of $R/I$, the claim follows.

From the claim we get that  $R/J$ and $R/I$ have the same finite
length so that the canonical surjection $R/J\to R/I$ is a bijection and $I=J.$

Finally, the given elements are a minimal basis of $I$ because the generators of $\n$ are analitically independent.

\medskip
\noindent{\bf We come now to the case $\tau(A)=h.$ }

\vskip 2mm
\noindent
In the case the Cohen-Macaulay type of $A$ is $h$, the maximum allowed,
we get by (\ref{m}) $\m=(y_1,y_2,\dots,y_h)$ where $(y_2,\dots,y_h)\subseteq 0:\m.$
This implies that  $y_1y_i=0$ for every $i=2,\dots,h$ and $y_iy_j=0$ for every $2\le i\le j\le h.$
Further we also have $y_1^{s+1}=0.$ The conclusion follows as in case i), but is
 even easier because the generators of $J$ are monomials.
 \end{proof}

 \begin{remark}
 It is clear that, for a stretched  local ring $A=R/I$ of maximal type,
 the minimal set of generators of $I$ we have found in the above theorem
 are a standard basis for $I.$ Namely we have that $I^*$ is the ideal
 generated by $X_1^{s+1}$ and the degree two monomials in $S,$ except for $X_1^2.$
This is not true in the case $\tau(A)<h.$
In this case,   the initial forms of the generators of $I$ in
$S=gr_{\n}(R)=\res[X_1,X_2,\dots,X_h]$ are the degree two monomials in $S,$
except for $X_1^2.$ The ideal $I^*$ is, as before,
the ideal generated by $X_1^{s+1}$ and the degree two monomials in $S,$ except for $X_1^2.$
 \end{remark}

 \begin{remark}
 It is clear that, given two  integers $1\le \tau\le h$ and a  regular local ring
 $(R,\n)$ with maximal ideal $\n$ minimally generated by $(x_1,x_2,\dots,x_h),$
 the ideals $I$ generated as in Theorem \ref{first} have the property that
 $A:=R/I$ is a stretched local ring of type $\tau.$
 \end{remark}

\vskip 3mm We have proved that if $R/I$ is a stretched Artinian
local ring of embedding dimension $h$, Cohen-Macaulay type $\tau<h$ and
socle degree $s,$ then we can find a minimal system of generators $x_1,\dots, x_h$ of $\n$ such that
$$I=(\{x_ix_j\}_{1\le i<j\le h},
\{x_j^2\}_{2\le j \le \tau},
\{x_i^2-u_ix_1^s\}_{ \tau +1\le i \le h})
$$ where the $u_i$ are units in $R.$
For every $\underline{u}=(u_j)_{j=\tau +1,\dots,h},$ we let $I(\underline{u})$ such an ideal.

\vskip 2mm We will use several time the following easy and well known Lemma that is a
consequence of Hensel's Lemma.

\begin{lemma}\label{crucial}
 Let $(A,\m)$ be an Artinian local ring with residue field $\res$ and let $a$
 be  an element in $A$ such that $\overline{a}\in \res^*.$
 If $\overline{b}^n=\overline{a}$ for some $\overline{b}\in \res,$
 then $c^n=a$ for some $c\in A, c\notin \m.$
\end{lemma}

\begin{proposition}
Let $I(\underline{u})$ as before and  assume that the residue
field $\res=R/\n$ verifies $\res^{1/2}\subseteq \res.$
 Then we can find a system of generators $y_1,\dots,y_h$ of $\n$
 such that
 $$
 I(\underline{u})=(\{y_iy_j\}_{1\le i<j\le h},
 \{y_j^2\}_{2\le j \le \tau},
 \{y_i^2-y_1^s\}_{ \tau +1\le i \le h}).
 $$
 \end{proposition}
\begin{proof}
Since $\res^{1/2}\subseteq \res,$ by the above Lemma we can find,
for every $i=\tau +1,\dots,h,$ elements $v_i \in R$ such that
$v_i^2\cong 1/u_i$ mod $I(\underline{u}).$
Hence $v_i\notin \n$ and we get
$$
v_i^2x_i^2-x_1^s\cong (1/u_i)x_i^2-x_1^s=(1/u_i)(x_i^2-u_ix_1^s)\cong 0.
$$
This proves that if we let
$$y_i=x_i, \ \  \text {for}\ \  i=1,\dots,\tau, \ \ \ y_i=v_ix_i\ \  \text{for}\ \  i=\tau +1,\dots,h,$$
then
$$(\{y_iy_j\}_{1\le i<j\le h},\{y_j^2\}_{2\le j \le \tau},\{y_i^2-x_1^s\}_{ \tau +1\le i \le h})
\subseteq
I(\underline{u}).
$$
Since the two ideals have the same Hilbert function, they must coincide.
\end{proof}


\bigskip
\section{Almost stretched Gorenstein  local rings}
In this section we are considering Artinian local rings  $(A,\m)$ such that
 the square of the maximal ideal  is minimally  generated by two elements.
 Recall that in Section 1 such a ring $A$ has been called almost stretched.
 If $A$ is almost stretched and Gorenstein, the Hilbert function
  of $A$ is given by

\begin{center}
\begin{tabular}{|c|c|c|c|c|c|c|c|c|c|c|c|}
0 & 1  & 2 & \dots & t & t+1  &\dots & s  & s+1  \\ \hline
1  & h  & 2 & \dots & 2 & 1 &  \dots & 1  &  0  \\
\end{tabular}
\end{center}

\noindent
with $h\ge 2$ and $s\ge t+1\ge 3.$

The structure result for almost stretched Gorenstein local rings will be a consequence of the following theorem.

\begin {theorem}\label{main}  Let $(A,\m)$ be an Artinian local ring which is Gorenstein with embedding dimension $h.$
If $A$ is almost stretched, then we can find integers $s\ge t+1\ge 3$ and a minimal basis $x_1,\dots,x_h$
of $\m$ such that
$$
\begin{cases}
x_1x_j=0        & \text {for $j=3,\dots,h$}\\
x_ix_j=0        &  \text{for $2\le i < j \le h$}\\
x_j^2=u_jx_1^s  &  \text{for $j=3,\dots,h$}\\
 x_2^2=ax_1x_2+wx_1^{s-t+1}\\
 x_1^tx_2=0
 \end{cases}
 $$
with suitable $w,u_3,\dots,u_h\notin \m$ and $a\in A.$
\end{theorem}
\begin{proof}
By Proposition \ref{lean} we may assume that  $\m=(x_1,\dots,x_h)$ with $$\m^j=\begin{cases}
   (x_1^j,x_1^{j-1}x_2) & j=2,\dots,t \\
     (x_1^j)  & j=t+1,\dots,s.
\end{cases}$$
We claim that we may assume also $(x_3,\dots,x_h)\subseteq (0):x_1.$
Namely, for $j\ge 3,$ we can write $x_1x_j=b_jx_1^2+c_jx_1x_2,$ hence
$x_1(x_j-b_jx_1-c_jx_2)=0.$ We get the claim by replacing   $x_j$ with $x_j-b_jx_1-c_jx_2$
for every $j\ge 3$.
This means that we have \begin{equation} \label{x1xj} {x_1x_3=x_1x_4=\dots=x_1x_h=0}.\end{equation}

\noindent
Further, since $\m^{t+1}=(x_1^{t+1}),$  for some $c\in A$, we have
\begin{equation}
\label{x-1^tx-2-}
x_1^tx_2=c x_1^{t+1}.
\end{equation}

\noindent
Let $y_2:=x_2-cx_1,$ then $$x_1^ty_2=x_1^t(x_2-cx_1)=x_1^tx_2-cx_1^{t+1}=0.$$  Since  $x_2$
 is not involved in equations (\ref{x1xj}), we may replace $x_2$ with $y_2$ in the generating
 set of $\m.$ Hence   we may assume that
\begin{equation}
\label{x-1^tx-2}
 { x_1^tx_2=0}.
\end{equation}
We notice that $  x_1^{t-1}x_2\notin \m^s$, otherwise $x_1^{t-1}x_2\in\m^{t+1},$ a
contradiction to the fact that $x_1^{t-1}x_2,x_1^t$ is a minimal basis of $\m^t.$
This implies that $x_1^{t-1}x_2$ cannot be in the socle of $A.$
Since by (\ref{x-1^tx-2}) and (\ref{x1xj})
\begin{equation*}
x_1^{t-1}x_2\in (0):(x_1,x_3,\dots,x_h),
\end{equation*}
we must have
 \begin{equation}\label{Gor}
 x_1^{t-1}x_2^2 \not=0.
 \end{equation}

 We want to prove now that we can find $a\in A$, $w\notin \m$ such that
 $$
 x_2^2=ax_1x_2+wx_1^{s-t+1}.
 $$
 In order to prove this we need the following easy remarks.

 \vskip 2mm \noindent
 {\bf Claim} 1. If for some $r,p \in A$ and $n\ge 2$ we have $x_2^2=rx_1x_2+px_1^n, $ then $n\le s-t+1.$
 If further $p\not\in \m,$ then $n= s-t+1.$

 \vskip 2mm \noindent
 Proof of Claim 1.
 We have
 $$
 x_1^{t-1}x_2^2=x_1^{t-1}(rx_1x_2+px_1^n)=px_1^{n+t-1}
 $$
 because by (\ref{x-1^tx-2}) $x_1^tx_2=0.$
 Since by (\ref {Gor}) $x_1^{t-1}x_2^2 \not=0,$ this implies $n+t-1\le s.$
 We have also $px_1^n=x_2(x_2-rx_1),$ hence,
 if $p\notin \m,$  $x_1^n=vx_2$ for some $v\in A.$
 As a consequence we get $x_1^{n+t}=vx_1^tx_2=0.$
 Since $x_1^s\not= 0,$ we have $n+t\ge s+1$ and the conclusion follows.

\vskip 3mm
\noindent {\bf Claim} 2.  If for some  $n\ge 2,$
$a\in A$ and $b\in \m,$  we have $x_2^2=ax_1x_2+bx_1^n$
then for some $c,d \in A$ we have $x_2^2=cx_1x_2+dx_1^{n+1}.$

\vskip 2mm \noindent
Proof of Claim 2 . This is easy because by (\ref{x1xj})
$x_1x_j =0$ for every $j\ge 3.$

\vskip 3mm
\noindent {\bf Claim} 3.
If for some   $a,b\in A$   we have  $x_2^2=ax_1x_2+bx_1^{s-t+1}$  then $b\notin \m.$

\vskip 2mm \noindent
Proof of Claim 3.  If, by contradiction, $b\in \m,$ then by Claim 2 and 1  we get
$$
s-t+2\le s-t+1.
$$
Since $\m^2=(x_1^2,x_1x_2),$   we have $ x_2^2=ax_1x_2+b x_1^2$ for some $a,b \in A$.
Thus,   as a trivial consequence of these three  claims, we get that for some $ a\in A$ and $w\notin \m$
\begin{equation}
\label{quadr}
x_2^2=ax_1x_2+wx_1^{s-t+1}.
\end{equation}

Now we recall that for every $j\ge 3,$ we have by (\ref{x1xj})
$$
x_j \m^2=x_j(x_1^2,x_1x_2)=0,
$$
so that, by using the Gorenstein assumption, we get
\begin{equation}
\label{nuova}
x_j \m \subseteq (0):\m=(x_1^s).
\end{equation}

Let us consider the ideal $J:=(x_3,\dots,x_h).$
By (\ref{nuova}), for every $3\le i \le j \le h,$
we have $x_ix_j= u_{ij}x_1^s$ with $u_{ij}\in A.$
We notice that if we have also $x_ix_j= w_{ij}x_1^s$,
then $(u_{ij}-w_{ij})x_1^s= 0,$ which implies $u_{ij}-w_{ij}\in \m.$

Hence we may define an inner product in  the  $\res=A/\m$-vector
space $V:=J/J\m$ by letting
$$
<\overline{x_i},\overline{x_j}>:=\overline{u_{ij}}\in A/\m
$$ and extending this definition  by bilinearity to $V\times V.$

Since the characteristic of $\res$ is not two,
the inner product can be diagonalized.
This means that we can find minimal  generators
$y_3,\dots,y_h$ of $J$ such that $y_iy_j=0$ for $i\not= j.$ If we replace
$x_3,\dots,x_h$ with $y_3,\dots,y_h$ in the generating set of $\m,$
it is clear that equations (\ref{x1xj}), (\ref{x-1^tx-2}), (\ref{quadr}) and (\ref{nuova})
are still valid.
Thus generators $x_1,\dots,x_h$ of $\m$ can be chosen so that
\begin{equation}
\label{aij=0}
x_ix_j=0
\end{equation}
for every $i$ and $j$ such that $3\le i < j \le h.$

\vskip 3mm
From (\ref{nuova}) and  for every $j\ge 3$ we have  $$x_j^2= u_jx_1^s$$ with $u_j\in A.$
We claim  that  $ u_j\notin \m$ for every $j\ge 3.$

In order to prove this claim, let us remember that again by
(\ref{nuova}) we have
$$x_2x_j= a_jx_1^s$$
for every $j\ge 3$ and suitable $a_j\in A.$
We fix $j\ge 3$ and let $$\rho :=wx_j-a_{j}x_1^{t-1}x_2.$$
Since $w\notin \m,$ it is clear that $\rho\notin \m^2$ so that
$\rho\notin \m^s\subseteq \m^2.$
This implies that  $\rho$ {\bf cannot be in the socle of $A$}.
We will use the following equalities:

\begin{center}
\begin{tabular}{lll}
  $x_1x_j= 0$ & for $j\ge 3$& see (\ref{x1xj}) \\ \\
  $x_1^tx_2= 0$ & &see (\ref{x-1^tx-2})             \\ \\
  $x_2^2= ax_1x_2+wx_1^{s-t+1}$&  &see (\ref{quadr}) \\ \\
  $x_jx_k= 0$&  for $3\le j<k\le h$ &see (\ref{aij=0})\\ \\
\end{tabular}
\end{center}

\noindent
We have

\medskip
\noindent
$\rho x_1=wx_1x_j-a_{j}x_1^tx_2=0,$

\medskip
\noindent
$\rho x_2=wx_2x_j-a_{j}x_1^{t-1}x_2^2= w a_{j}x_1^s-a_{j}x_1^{t-1}(ax_1x_2+wx_1^{s-t+1})= wa_{j}x_1^s-wa_{j}x_1^s= 0,$

\medskip
\noindent
$\rho x_k=wx_jx_k-a_{j}x_1^{t-1}x_2x_k=0$ \quad if $k\ge 3,$ $k\not= j,$

\medskip
\noindent
$\rho x_j=wx_j^2-a_{j}x_1^{t-1}x_2x_j= wu_jx_1^s.$

\medskip
\noindent
Since  $\rho$  cannot be in the socle,  we must have $ u_j\notin \m.$ This proves the Claim.

\medskip
As a consequence we may assume that for every $j\ge 3$ and suitable $ u_j\notin \m$ we have
\begin{equation}
\label{xjxj}
  x_j^2=u_jx_1^s .
 \end{equation}

We come now to the last manipulation of our elements.
As a consequence of the above claim,
we may consider the element
$$
y_2:=x_2-\sum_{i=3}^h u_i^{-1}a_{i}x_i.
$$
For every $j=3,\dots,h$ we have by using (\ref{aij=0})
$$
y_2x_j=x_2x_j-\sum_{i=3}^h u_i^{-1}a_{i}x_ix_j= a_{j}x_1^s-u_j^{-1}a_{j}x_j^2= a_{j}x_1^s-u_j^{-1}a_{j}u_jx_1^s= 0.
$$
Further we have
$$
x_1^tx_2=x_1^t(y_2+\sum_{i=3}^h u_i^{-1}a_{i}x_i)= x_1^ty_2.
$$
Finally  let $d:=x_2-y_2=\sum_{i=3}^h u_i^{-1}a_{i}x_i.$
Then $d \in J:=(x_3,\dots,x_h)$ and we have
$$x_1d= 0 \ \ \ \ \ \  y_2d= 0.$$
Since, by (\ref{nuova}), $J\m\subseteq (x_1^s),$
we have
$$d^2 = px_1^s$$
for some $p\in A.$
It follows that
$$x_2^2-ax_1x_2-wx_1^{s-t+1}=(y_2+d)^2-ax_1(y_2+d)-wx_1^{s-t+1}= y_2^2+d^2-ax_1y_2-wx_1^{s-t+1}=$$
$$=  y_2^2-ax_1y_2-wx_1^{s-t+1}+px_1^s= y_2^2-ax_1y_2-(w-px_1^{t-1})x_1^{s-t+1}$$
where $w-px_1^{t-1}\notin \m.$

Thus we may replace $x_2$ with $y_2$ and  finally we get  a basis
$x_1,\dots,x_h$ for $\m$ so that

$$ \begin{cases}
x_1x_j=0    & \text {for $j=3,\dots,h$}\\
x_ix_j=0     &  \text{for $2\le i < j \le h$}\\
x_j^2=u_jx_1^s &  \text{for $j=3,\dots,h$}\\
 x_2^2=ax_1x_2+wx_1^{s-t+1}\\
 x_1^tx_2=0
 \end{cases} $$
with suitable $w,u_3,\dots,u_h\notin \m$ and $a\in A.$
\end{proof}

\vskip 5mm
As a consequence of this theorem we get a structure theorem for almost stretched Artinian and Gorenstein local rings.

\begin{corollary}
\label{basic}
Let $(R,\n)$ be a regular local ring of dimension $h$ and $I\subseteq \n^2$
 an ideal such that $(A=R/I,\m=\n/I)$ is almost stretched Artinian and Gorenstein.
 Then there is a minimal basis $x_1,\dots,x_h$ of $\n$ such that $I$ is
 minimally generated by the elements
 $$
 \{x_1x_j\}_{j=3,\dots,h}\ \ \{x_ix_j\}_{2\le i < j \le h} \ \  \{x_j^2-u_jx_1^s\}_{j=3,\dots,h}\ \
  x_2^2-ax_1x_2-wx_1^{s-t+1}, \ \ \ \ x_1^tx_2.
$$
with $w,u_3,\dots,u_h \notin \n$ and $a\in R.$
\end{corollary}
\begin{proof} By the \thmref{main}  we can find a basis $x_1,\dots,x_h$ of $\n$ such that the ideal
$J$ generated by the above elements is contained in $I$.
We need  to show that $I$ is indeed equal to $J.$
We first remark that modulo $J$ we have
$$x_1^{s+1}=x_1^tx_1^{s-t+1}\cong x_1^t\frac{x_2^2-ax_1x_2}{w} \cong
x_1^tx_2\frac{x_2-ax_1}{w}\cong 0 $$
so that $x_1^{s+1}\in J.$

Passing to the ideals of initial forms in the polynomial ring
$$S=gr_{\n}(R)=\oplus_{j\ge 0}(\n^j/\n^{j+1})=(R/\n)[X_1,\dots,X_h],$$
we have $$I^*\supseteq J^*\supseteq K$$ where $K$ is the ideal in $S$
generated by the elements
$$ \{X_1X_j\}_{j=3,\dots,h}\ \   \{X_iX_j\}_{2\le i < j \le h} \ \
 \{X_j^2\}_{j=3,\dots,h}\ \ X_1^tX_2, \ \ \ X_1^{s+1}$$
 and the quadric $Q:= X_2^2-\overline{a}X_1X_2$ in the case  $s\ge t+2,$
 or  $Q:=X_2^2-\overline{a}X_1X_2-\overline{w}X_1^2$  in the case $s=t+1.$

In both cases we have $X_j S_1\subseteq K$ for every $j\ge 3$ so
that $$(K+(X_3,\dots, X_h))_n=K_n$$ for every $n\not= 1.$
This implies that for every $n\not= 1$
$$
H_{S/K}(n)= H_{S/(K+(X_3,\dots, X_h))}(n)=H_{\res[X_1,X_2]/(Q,X_1^tX_2,X_1^{s+1})}(n).
$$

Now we compute the Hilbert Function of the graded algebra $\res[X_1,X_2]/(Q,X_1^tX_2,X_1^{s+1}).$
We let $B:=\res[X_1,X_2];$ in the case $Q= X_2^2-\overline{a}X_1X_2=X_2(X_2-\overline{a}X_1)$,
we have an exact sequence
$$
0\to B/(X_2-\overline{a}X_1,X_1^t)(-1)   \overset{X_2}{\to}B/(Q,X_1^tX_2)\to B/(X_2)\to 0
$$
which enables us to compute the Hilbert Series of  $B/(Q,X_1^tX_2)$:
$$P_{B/(Q,X_1^tX_2)}(z)=zP_{B/(X_2-\overline{a}X_1,X_1^t)}(z)+P_{B/(X_2)}(z)=$$
$$=\frac{z(1-z)(1-z^t)+(1-z)}{(1-z)^2}=\frac{1+z-z^{t+1}}{1-z}$$ which gives the Hilbert Function

\begin{center}
\begin{tabular}{|c|c|c|c|c|c|c|c|c|c|c|c|}
0 & 1  & 2 & \dots & t & t+1 & \dots & s  & s+1 &s+2&  \dots \\
\hline 1  & 2  & 2 & \dots & 2 & 1&  \dots & 1  &  1 &1 &\dots \\
\end{tabular}
\end{center}

Since $X_1^{s+1}\notin (Q,X_1^tX_2),$ the Hilbert Function of $\res[X_1,X_2]/(Q,X_1^tX_2,X_1^{s+1})$ is
\begin{center}
 \begin{tabular}{|c|c|c|c|c|c|c|c|c|c|c|c|c|}
  0 & 1  & 2 & \dots & t & t+1 & t+2 &\dots & s  & s+1  \\ \hline
  1  & 2  & 2 & \dots & 2 & 1&1 &  \dots & 1  &  0  \\
  \end{tabular}
  \end{center}
   so that the Hilbert Function of $S/K$ is
\begin{center}
\begin{tabular}{|c|c|c|c|c|c|c|c|c|c|c|c|c|}
0 & 1  & 2 & \dots & t & t+1 & t+2 &\dots & s  & s+1  \\ \hline
1  & h  & 2 & \dots & 2 & 1&1 &  \dots & 1  &  0  \\
\end{tabular}
\end{center}
the same as that of $S/I^*.$

In the case $s=t+1$ we have $Q=X_2^2-\overline{a}X_1X_2-\overline{w}X_1^2$
with $\overline{w}\not= 0.$
Hence  $\{Q,X_1^tX_2\}$ is a regular sequence and $\res[X_1,X_2]/(Q,X_1^tX_2)$ has Hilbert Function

\begin{center}
\begin{tabular}{|c|c|c|c|c|c|c|c|c|c|c|c|}
0 & 1  & 2 & \dots & t & t+1=s & t+2   \\ \hline
1  & 2  & 2 & \dots & 2 & 1&    0  \\
\end{tabular}
\end{center}

We remark that in this case we have $X_1^2\in (Q,X_2)$ so that
$$
X_1^{s+1}=X_1^{t+2}=X_1^tX_1^2\in (Q,X_1^tX_2).
$$

In any case we have proven that $S/I^*$ and $S/K$ have the same Hilbert Function.
This implies that $I^*=J^*=K$ so that  the Hilbert Function of $R/I$ and $R/J$ are the same.
Hence  $R/I$ and $R/J$ have the same finite length, so the canonical epimorphism $R/J \to R/I$
is an isomorphism and $I=J$ as claimed.
\end{proof}

\begin{remark}
Notice that in the proof of \corref{basic} we describe the ideal $I^*$: is generated by
$$ \{X_1X_j\}_{j=3,\dots,h}\ \   \{X_iX_j\}_{2\le i < j \le h} \ \
 \{X_j^2\}_{j=3,\dots,h}\ \ X_1^tX_2, \ \ \ X_1^{s+1}$$
 and the quadric $Q:= X_2^2-\overline{a}X_1X_2$ in the case  $s\ge t+2,$
 or  $Q:=X_2^2-\overline{a}X_1X_2-\overline{w}X_1^2$  in the case $s=t+1$,
 with $ \overline{w}\neq 0 , \overline{a}\in\res$.
\end{remark}

\bigskip
We want to prove now the converse of the above result.
Notice that for the next Lemma we even do not need neither regular nor local.

\begin{lemma}
\label{primo}
Let $B$ a ring, $t\ge 2$, $h\ge 2,$ $s\ge t+1$ and  $\n=(x_1,\dots,x_h)$ an ideal in $B$.
Let  $J$ be the ideal generated by
$$\{x_1x_j\}_{j=3,\dots,h}\ \ \{x_ix_j\}_{2\le i < j \le h} \ \
\{x_j^2-u_jx_1^s\}_{j=3,\dots,h}\ \  x_2^2-ax_1x_2-wx_1^{s-t+1}, \ \ \ \ x_1^tx_2.
$$
If $w$ is a unit in $B,$  then $$\n^{s+1}\subseteq J.$$
\end{lemma}
\begin{proof}  For every $i\not= j,$ save  for $(i,j)=(1, 2),$ we have $$ x_ix_j\in J.$$
For every $3\le j \le h,$ $$x_j^2\in J+(x_1^s),$$  and since $s-t+1\ge 2,$ $$x_2^2\in J+(x_1^2,x_1x_2).$$

We claim that for every $r\ge 2$ we have $$\n^r\subseteq J+(x_1^r,x_1^{r-1}x_2).$$

If $r=2,$ we have $\n^2\subseteq J+(x_1^2,x_1x_2)$ by the above three properties.
Let us proceed by induction on $r$. We have $$\n^{r+1}=\n \n^{r}\subseteq J+\n (x_1^r,x_1^{r-1}x_2)=$$
$$=J+(x_1,x_2)(x_1^r,x_1^{r-1}x_2)=J+(x_1^{r+1},x_1^rx_2,x_1^{r-1}x_2^2).$$
The claim follows because $x_2^2\in J+(x_1^2,x_1x_2)$, so
$$
x_1^{r-1}x_2^2\in J+(x_1^{r+1},x_1^rx_2).
$$
From the claim we have $\n^{s+1}\subseteq J+(x_1^{s+1},x_1^s x_2).$
Since $s\ge t,$ we get $x_1^sx_2\in(x_1^tx_2)\subseteq J;$  on the other hand,
since $w$ is a unit   we get modulo $J$ the equalities
$$
x_1^{s+1}=(x_1^t/w)wx_1^{s-t+1} \cong (x_1^t/w)(x_2^2-ax_1x_2)\cong 0 .
$$
The conclusion follows.
\end{proof}

\medskip
We come now to a very crucial step in our way.

\begin{lemma}
\label{secondo} Let $R$ be a regular local ring of dimension $h\ge
2,$ $\n=(x_1,\dots,x_h)$ the maximal ideal of $R,$ $s\ge t+1\ge 3$
and $a,u_3,\dots, u_h,w \in R.$ Let $I$ be the ideal generated by
$$
\{x_1x_j\}_{j=3,\dots,h}\ \ \{x_ix_j\}_{2\le i < j \le h} \ \  \{x_j^2-u_jx_1^s\}_{j=3,\dots,h}\ \
q:=x_2^2-ax_1x_2-wx_1^{s-t+1},  x_1^tx_2.
$$
If $u_3,\dots, u_h,w \notin \n,$   then

\vskip 2mm $(1)$ $\overline{x_1}^t,  \overline{x_1}^{t-1} \overline{x_2}\in (\n^t+I)/(\n^{t+1}+I)$ are
$(R/\n)$-linearly independent elements,

\vskip 2mm
$(2)$ $x_1^s\notin I.$
\end{lemma}
\begin{proof}
In order to prove $(1)$ we need to show that  if
$\lambda x_1^t+\mu x_1^{t-1}x_2\in I+\n^{t+1},$
then $\lambda, \mu \in \n.$
It is clear that if $\lambda x_1^t+\mu x_1^{t-1}x_2\in I+\n^{t+1},$ then
$$
\lambda x_1^t+\mu x_1^{t-1}x_2\in I+\n^{t+1}+(x_3,\dots,x_h)=
(x_3,\dots,x_h)+(x_1,x_2)^{t+1}+(x_1^s,x_1^tx_2,q)
$$
$$=(x_3,\dots,x_h)+(x_1,x_2)^{t+1}+(q).
$$
Let's read the above condition in the two dimensional regular local ring $T:=R/(x_3,\dots,x_h)$,
whose maximal ideal is generated by the residue class of $x_1$ and $x_2$ modulo $(x_3,\dots,x_h).$
By abuse of notation, we again denote these elements by $x_1,x_2$  and the maximal ideal of $T$ by $\n.$
We have $$\lambda x_1^t+\mu x_1^{t-1}x_2=eq+z$$ where $z\in \n^{t+1}.$
This implies that $eq\in \n^t.$
If   $eq\in \n^{t+1},$   the conclusion follows by the analytic independence of $x_1,x_2.$
If $eq\notin  \n^{t+1},$ then  since $q=x_2^2-ax_1x_2-wx_1^{s-t+1}\in \n^2,$
we have $e\in \n^{t-2},$ $e\notin \n^{t-1}.$
By passing to the associated graded ring $(T/\n)[X_1,X_2]$ of $T,$ we get
$$
X_1^{t-1}(\overline{\lambda} X_1+\overline{\mu} X_2)=e^* q^*.
$$
Since  $X_1$ is not a factor of  $q^*$, $X_1^{t-1}$ must be a factor of $e^*.$
This is a contradiction because $e^*$ is an homogeneous element of degree $t-2.$
 The conclusion follows.

\medskip
Let us prove $(2)$.
 By contradiction, let
 $$
 x_1^s=\sum_{j=3}^h\lambda_jx_1x_j+ \sum_{j=3}^h\rho_j(x_j^2-u_jx_1^s)+
 \sum_{2\le i<j\le h}\mu_{ij}x_ix_j+\sigma x_1^tx_2+\alpha q.
 $$
 Since $s\ge t+1\ge 3,$  this implies
 $$
 \sum_{j=3}^h\lambda_jx_1x_j+ \sum_{j=3}^h\rho_j x_j^2+\sum_{2\le i<j\le h}\mu_{ij}x_ix_j+
 \alpha(x_2^2-ax_1x_2-wx_1^{s-t+1})\in \n^3.
 $$
 By the analytic independence of $x_1,\dots,x_h$, all the coefficients of the degree two monomials
in $x_1,\dots,x_h$ must be in $\n.$
In particular $\rho_j \in \n$ for every $j=1,\dots,h.$
This implies that
$$
x_1^s\in (x_3,\dots,x_h)+(x_1^tx_2,q)+\n^{s+1}.
$$
As we did before, we pass to  the two dimensional regular local ring $T:=R/(x_3,\dots,x_h)$
 whose maximal ideal is still denoted by $\n$ and generated by $x_1,x_2.$
 We can write
 \begin{equation}
\label{vvv} x_1^s=\sigma x_1^tx_2+\alpha q+\beta
\end{equation}  where $\beta\in \n^{s+1}.$
This implies that $x_1^s+\alpha wx_1^{s-t+1}\in (x_2,x_1^{s+1})$ so
that  we can write $x_1^s+\alpha wx_1^{s-t+1}=x_2a+x_1^{s+1}b$ for some
$a,b\in T.$ This gives $$x_1^{s-t+1}(x_1^{t-1}+\alpha w-bx_1^t)=x_2a.$$
Since $x_1^{s-t+1}, x_2$ is a regular sequence in $T,$
we get $x_1^{t-1}+\alpha w-bx_1^t=x_2c$ for some $c\in T.$
Hence $\alpha w=x_1^{t-1}(bx_1-1)+x_2c$ and since $w$ is a unit,
we finally get $$\alpha=vx_1^{t-1}+dx_2$$ for some $v,d \in T,$ $v\notin \n.$
Let us use this formula in equation (\ref{vvv}). We get
\begin{equation}
\label{www} x_1^s=\sigma x_1^tx_2+ (vx_1^{t-1}+dx_2)q+\beta
\end{equation}
where $\beta\in \n^{s+1}$ and $v\notin \n.$

We claim now that if for some $r\ge 2$ and $j\ge 2$ we have, as in (\ref{www})
with $j=s$ and $r=t$,
$$
x_1^j-\sigma x_1^rx_2- (vx_1^{r-1}+dx_2)q\in \n^{j+1}
$$
then, for suitable  $e\in T,$  we get also
  $$
  x_1^{j-1}-\sigma x_1^{r-1}x_2- (vx_1^{r-2}+ex_2)q\in \n^j.
  $$
  Since $q=x_2^2-ax_1x_2-wx_1^{s-t+1},$
  the assumption of the claim implies
  $$
  dx_2^3\in (x_1)+\n^{j+1}=(x_1)+(x_2^{j+1}).
  $$
  Now, since $j+1\ge 3$ and $x_1,x_2^3$ is a regular sequence,
  we get $d=ex_1+fx_2^{j-2}$ for some $e,f\in T$ so that
  $x_1^j-\sigma x_1^rx_2- (vx_1^{r-1}+ex_1x_2)q\in \n^{j+1}.$
  Since $\n^{j+1}\cap (x_1)=x_1\n^j,$ it follows that
  $$
  x_1^{j-1}-\sigma x_1^{r-1}x_2- (vx_1^{r-2}+ex_2)q\in \n^j
  $$
  and the claim is proved.

Starting from (\ref{www}), where we let  $j=s$ and $r=t,$
we apply $t-1$ times the claim and we get
$$
x_1^{s-t+1}-\sigma x_1x_2- (v+gx_2)q\in \n^{s-t+2}
$$
for some $g\in T.$
This implies
$$
(v+gx_2)x_2^2\in (x_1)+\n^{s-t+2}=(x_1,x_2^{s-t+2}),
$$ so that, since $s-t+2\ge 3,$ we get  $vx_2^2\in (x_1,x_2^3),$
which is a contradiction because $v\notin \n.$
\end{proof}

\begin{corollary}
\label{HF}
Let $R$ be a regular local ring of dimension $h\ge 2,$
$\n=(x_1,\dots,x_h)$ the maximal ideal of $R,$
$s\ge t+1\ge 3$ and $a,u_3,\dots, u_h,w \in R.$
Let $I$ be the ideal generated by
$$
\{x_1x_j\}_{j=3,\dots,h}\ \ \{x_ix_j\}_{2\le i < j \le h} \ \  \{x_j^2-u_jx_1^s\}_{j=3,\dots,h}\ \
q:=x_2^2-ax_1x_2-wx_1^{s-t+1},  x_1^tx_2.
$$
If $u_3,\dots, u_h,w \notin \n,$
then the Hilbert Function of $R/I$ is
\begin{center}
{\rm
\begin{tabular}{|c|c|c|c|c|c|c|c|c|c|c|c|c|}
0 & 1  & 2 & \dots & t & t+1 & t+2 &\dots & s  & s+1  \\ \hline
1  & h  & 2 & \dots & 2 & 1&1 &  \dots & 1  &  0  \\
\end{tabular}
}
\end{center}
\end{corollary}
\begin{proof}
We have seen in the proof of  \lemref{primo} that $\n^r\subseteq
J+(x_1^r,x_1^{r-1}x_2)$ for every $r\ge 2.$ This proves that all
the powers of $\n/I$ can be generated by two elements. By a) of
Lemma \ref{secondo} we get $H_{R/I}(t)=2$, which implies, by the
characterization of Hilbert functions due to Macaulay,
$H_{R/I}(j)=2$ for every $2\le j\le t.$ Since $x_1^tx_2\in I,$ we
also have $H_{R/I}(t+1)\le 1,$ which implies $H_{R/I}(j)\le 1$ for
every $j\ge t+1.$ The conclusion follows because $x_1^s \notin I$
and $\n^{s+1}\subseteq I.$
\end{proof}

\bigskip
We are ready to prove the converse of  \corref{basic}.

\begin{theorem}
\label{rbasic}
Let $R$ be a regular local ring of dimension $h\ge 2,$     $\n=(x_1,\dots,x_h)$
the maximal ideal of $R,$ $s\ge t+1\ge 3$ and $a,u_3,\dots, u_h,w \in R.$
Let $I$ be the ideal generated by
$$
\{x_1x_j\}_{j=3,\dots,h}\ \ \{x_ix_j\}_{2\le i < j \le h}
 \ \  \{x_j^2-u_jx_1^s\}_{j=3,\dots,h}\ \  x_2^2-ax_1x_2-wx_1^{s-t+1},  x_1^tx_2.
$$
If $u_3,\dots, u_h,w \notin \n,$
then $R/I$ is an almost stretched  Gorenstein local ring with Hilbert function
\begin{center}
{\rm
\begin{tabular}{|c|c|c|c|c|c|c|c|c|c|c|c|c|}
0 & 1  & 2 & \dots & t & t+1 & t+2 &\dots & s  & s+1  \\ \hline
1  & h  & 2 & \dots & 2 & 1&1 &  \dots & 1  &  0  \\
\end{tabular}
}
\end{center}
\end{theorem}
\begin{proof} After the above Corollary we need only to prove that $R/I$ is Gorenstein.

We let $\m:=\n/I$ and $y_i:=\overline{x_i}\in A=R/I.$
By \lemref{secondo} we have $\m^j=(y_1^j,y_1^{j-1}y_2)$ for every $j=2,\dots t,$
and $\m^j=(y_1^j)$ for $j=t+1,\dots,s.$
We prove the theorem in three steps.

\vskip 2mm
\noindent
{\bf Claim} 1. If for some $j\not= 1,t, s$ and some  $r\in \m^j$ we have $ry_1=0,$ then $r\in \m^{j+1}.$

\vskip 2mm \noindent
Proof of Claim 1. Let $2\le j \le t-1;$ then $r=\lambda y_1^j+\mu y_1^{j-1}y_2.$
We have  $$0=ry_1=\lambda y_1^{j+1}+\mu y_1^jy_2.$$
Since $y_1^{j+1}, y_1^jy_2$ is a minimal basis of $\m^{j+1},$ we have $\lambda,\mu \in \m$ and
$r\in \m^{j+1}.$ The case $t+1 \le j \le s-1$ is even easier.

\vskip 2mm
\noindent{\bf Claim} 2. If for some   $r\in \m^t$ we have $ry_1=ry_2=0,$ then $r\in \m^{t+1}.$

\vskip 2mm \noindent
Proof of Claim 2. Let $r=\lambda y_1^t+\mu y_1^{t-1}y_2.$
Since $y_1^ty_2=0,$ we have $0=ry_1=\lambda y_1^{t+1}.$
This implies $\lambda\in \m.$
On the other hand we have
$$
0=ry_2=\mu y_1^{t-1}y_2^2=\mu y_1^{t-1}(\overline{a}y_1y_2+\overline{w} y_1^{s-t+1})=\mu \overline{w}y_1^s.
$$
Since $\overline{w}$ is a unit in $A,$ this implies $0=\mu y_1^s$ so that $\mu \in \m.$ Thus $r\in \m^{t+1}.$

\vskip 2mm These two Claims prove that if $r\in \m^2$ and $ry_1=ry_2=0,$ then $r\in \m^s.$

\vskip 3mm
\noindent
{\bf Claim} 3. If $r\in (0):\m$ then $r\in \m^2,$ so that $r\in \m^s$ and  $A$ is Gorenstein.

\vskip 2mm \noindent
Proof of Claim 3. Let $r\in (0):\m$; then $r\in  \m $ and we can write  $r=\sum_{i=1}^h \lambda_iy_i.$
 Since $y_1y_j=0$ for every $j\ge 3,$ we have
 $$
 0=ry_1=\lambda_1y_1^2+\lambda_2y_1y_2.
 $$
 This implies $\lambda_1,\lambda_2 \in \m$ so that $r=\sum_{i=3}^h \lambda_iy_i+b$ with $b\in \m^2.$

Since  $y_2y_j=0$ for every $j\ge 3,$ we get  $0=ry_1=by_1$ $0=ry_2=by_2;$ by Claim 2
this implies $b\in \m^s.$ Since $y_iy_j=0$ for every
$3\le i< j\le h,$ and $\m^{s+1}=0,$ we get
$$
0=ry_j=\lambda_jy_j^2=\lambda_j\overline{u_j}y_1^s.
$$
Since $\overline{u_j}$ is a unit in $A,$ this implies   $\lambda_jy_1^s=0$ so that
$\lambda_j\in \m$ and   $r\in \m^2.$
The proof of the Claim 3 and of the theorem is complete.
\end{proof}

\bigskip

The structure theorem of almost stretched Gorenstein local rings we have proved,
can be refined under a mild assumption on the residue field of $R.$
This will be crucial  for the study of the moduli problem
 and it is a consequence of the main structure \thmref{main} and \lemref{crucial}.

\begin{proposition}
\label{model}
Let $(R,\n,k)$ be a regular local ring of dimension $h\ge 2,$  and $I$ an ideal
in $R$ such that $R/I$ is almost stretched Artinian and Gorenstein.
If  $\res^{1/2}\subseteq \res,$ then we can find integers
$s\ge t+1\ge 3,$  a minimal system of generators $x_1,\dots,x_h$ of $\n$
and an element $a\in R,$ such that $I$ is generated by
$$\{x_1x_j\}_{j=3,\dots,h}\ \ \{x_ix_j\}_{2\le i < j \le h} \ \  \{x_j^2-x_1^s\}_{j=3,\dots,h}\ \
x_2^2-ax_1x_2-x_1^{s-t+1},\ \   x_1^tx_2.
$$
\end{proposition}
\begin{proof}
We know that integers $s\ge t+1\ge 3$ can be found and a minimal system of generators
$y_1,\dots,y_h$ of $\n$  can be constructed in such a way that
$I$ is generated by
$$\{y_1y_j\}_{j=3,\dots,h}\ \ \{y_iy_j\}_{2\le i < j \le h} \ \  \{y_j^2-u_jy_1^s\}_{j=3,\dots,h}\ \
y_2^2-by_1y_2-wy_1^{s-t+1},  \ \ y_1^ty_2.
$$
with $w,u_3,\dots,u_h \notin \n$ and $b\in R.$
By \lemref{crucial} we can find elements $v,r_3,\dots, r_h\in R$ such that modulo $I$
we have
$$
v^2\cong (1/w),\ \  r_3^2\cong (1/u_3), \dots, r_h^2\cong (1/u_h).
$$
From this is clear that $v,r_3,\dots, r_h$ are units in $R$ and we can make the
following change of minimal generators for
$\n:$
$$
x_1=y_1,\ \  x_2=vy_2, \ \ x_3=r_3y_3,\ \ \dots,\ \ x_h=r_hy_h.
$$
We have
$$
y_2^2-by_1y_2-wy_1^{s-t+1}=(x_2^2/v^2)-bx_1(x_2/v)-wx_1^{s-t+1}\in I,
$$ hence
$x_2^2-bvx_1x_2-v^2wx_1^{s-t+1}\in I.$
Since $v^2w=1+d$ with $d\in I,$ if we let $a:=bv,$ we  get
$$
x_2^2-ax_1x_2-x_1^{s-t+1}\in I.
$$
Further for every $j=3,\dots,h$ we have
$$
y_j^2-u_jy_1^s=(x_j/r_j)^2-u_jx_1^s\in I,
$$
hence $x_j^2-r_j^2u_jx_1^s\in I.$
Since $r_j^2u_j=1+e$ with $e\in I,$ we get for every $j=3,\dots,h$
$$
x_j^2-x_1^s\in I.
$$

Hence $I$ contains the ideal generated by
$$
\{x_1x_j\}_{j=3,\dots,h}\ \ \{x_ix_j\}_{2\le i < j \le h} \ \
\{x_j^2-x_1^s\}_{j=3,\dots,h}\ \ x_2^2-ax_1x_2-x_1^{s-t+1},
x_1^tx_2.
$$
Since by  \corref{HF} these  two ideals have the same Hilbert function, they coincide.
\end{proof}


\bigskip
\section{Classification of Gorenstein local algebras with Hilbert function (1,2,2,2,1,1,1)}

\bigskip
We have seen in \secref{SLC} that the Cohen-Macaulay type determines
the moduli class of stretched Artinian local rings. In the case of
almost stretched Artinian local rings, the problem is not so easy,
even in the Gorenstein case. For example it has been proved in
\cite{CN07} that if $A$ is Gorenstein with Hilbert function
$1,2,2,1$, we have only two models, namely the ideals $I=(x^2,y^3)$
and $I=(xy,x^3-y^3).$ But already in the next  case with symmetric
Hilbert function $1,2,2,2,1,$ we have at least three different
models, namely two ideals which are homogeneous $I=(x^2,y^4)$,
$I=(xy,x^4-y^4)$ and one which is not homogeneous, the ideal
$I=(x^4+ 2 x^3 y, y^2-x^3)$.

But things become soon even more complicate, already in the complete intersection case,
the case $h=2.$
We are going to study the moduli problem for complete intersection
local rings with Hilbert function $1,2,2,2,1,1,1.$
We will see that in this case we have a one-dimensional family.

\bigskip
In the following, $(R,\n)$ is a two dimensional regular local ring such that $\res=R/\n$
has the property $ \res^{1/2} \subseteq \res ;$  $I$ is an ideal in $R$ such that $A=R/I$ is
Gorenstein with Hilbert function $1,2,2,2,1,1,1.$
We are not going into all the details, better we try simply to give an idea of what is going on.

By the main structure theorem we know that there exists a system of generators
$y_1,y_2$ of $\n$ and an element $a\in R$ such that, \propref{model},
$$
I=(y_1^3y_2,y_2^2-ay_1y_2-y_1^4).
$$

\vskip 3mm \noindent
{\bf Case 1: $a\notin \n.$ }
Let us change the generators as follows: $$z_1=ay_1-y_2,\ \ \ z_2=y_1^3+ay_2.$$
We have $$d:=det\ \ \begin{pmatrix}
     a & y_1^2   \\
   -1   &  a
\end{pmatrix}\ \ =a^2+y_1^2\notin \n$$ so that $z_1,z_2$ is a minimal system
of generators of $\n.$
We have $$z_1z_2=(ay_1-y_2)(y_1^3+ay_2)=-a(y_2^2-ay_1y_2-y_1^4)-y_1^3y_2\in I.$$
Since  $I$ contains the product of two minimal generators of $\n,$ then
there exists a system of generators $x,y$ of $\n$ such that
$$I=(xy,y^4-x^6).
$$

\vskip 3mm \noindent
{\bf Case 2: $a\in \n.$ }
In this case,  we write $a=by_1+cy_2,$ and choose $v\in R$ such that
$1-cy_1\cong v^2$ modulo $I$, \lemref{crucial}.
 Notice that $v\notin \n$,  so that we can change the generators as follows
 $$
 x_1=y_1, \ \ \ \ x_2=vy_2
 $$
 and prove that
 $$
 I=(x_1^3x_2,x_2^2-dx_1^2x_2-x_1^4)
 $$
 with $d=b v^{-1}\in R.$

\vskip 3mm \noindent
{\bf Case 2a: $d\in \n.$ }
In this case we write $d=f x_1+ex_2$ and choose $v\in R$ such that $v^2\cong 1-ex_1^2$ modulo $I.$
It is clear that $v\notin \n$ so that we can change the generators of $\n$
by letting
$$
x=x_1,\ \ \ \ y=vx_2.
$$
Then it is easy to    prove that
$$
 I=(x^3y,y^2-x^4).
$$

Let now consider the case $d\notin \n.$
We distinguish two subcases, $d^2+4\in \n$ and $d^2+4\notin \n$.
We first assume that

\vskip 3mm \noindent
{\bf Case 2b1: $d^2+4\in \n.$ }
In this case we have modulo $I$
$$
(x_2-(d/2)x_1^2)^2\cong x_1^4+(d^2/4)x_1^4\cong x_1^4(1+(d^2/4))=ex_1^5
$$
with $e\in R.$
It follows that if we let
$$
l:=x_2-(d/2)x_1^2+(e/d)x_1^3+(e^2/d^3)x_1^4
$$
then $l^2\in I.$
Modulo $I$ we have
$$
x_1^3l=x_1^3(x_2-(d/2)x_1^2+(e/d)x_1^3+(e^2/d^3)x_1^4)\cong -(d/2)x_1^5+(e/d)x_1^6=
$$
$$=x_1^5(-d/2+(e/d)x_1)=vx_1^5
$$
with $v\notin \n.$
It follows that $J=(l^2,x_1^3l-vx_1^5)\subseteq I$.
Next we  prove $J=I$.

Notice that $x, l$ form a minimal system of generators of $\n$ and we denote by $L$ the initial form
of $l$ in the associated graded ring $gr_{\n}(R)$.
In order to prove that $I=J$  we need to show that the Hilbert function of $R/J$ is $1,2,2,2,1,1,1.$
We have
$$
(X^3L,L^2)\subseteq J^*\subseteq I^*,
$$
so we have to prove that
$$
Z^7\in J.
$$
Notice that, modulo $J$, we have
$$
v x^7 = x^5 l = v^{-1} (x^3 l^2)=0.
$$
Hence $x^7\in J$, so $(l^2,x_1^3l-vx_1^5)=I.$

Now we have
$$
(l^2,x_1^3l-vx_1^5)=(l^2,(x_1^3l/v)-x_1^5)=((l/v)^2,x_1^3(l/v)-x_1^5).
$$
If we let $x:=l/v, y=x_1$ then $\n=(x,y)$ and
$$
 I=(x^2,xy^3-y^5).
$$

\vskip 3mm \noindent
{\bf Case 2b2: $d^2+4\notin \n.$ }
We can find $c,e\in R\setminus \n$ such that modulo $I$ we have $ c^2  \cong d^2+4$ and $e^2\cong -(2/c)$,
\lemref{crucial}.
We let $p:=d/c$ and change the generators of $\n$ by letting
$$
x=(x_1/e),\ \ \ \ y=x_2+p(x_1/e)^2.
$$
We get
$$
x_1=xe, \ \ \ \ \   x_2=y-px^2
$$
so that modulo $I$ we get
$$
 0\cong x_1^3x_2=x^3e^3(y-px^2)=e^3(x^3y-px^5)
 $$
which implies $x^3y-px^5\in I.$
Further
$$
0\cong x_2^2-dx_1^2x_2-x_1^4= (y-px^2)^2-dx^2e^2(y-px^2)-x^4e^4=
$$
$$
=y^2-x^2y(2p+de^2)+x^4(p^2+de^2p-e^4)\cong y^2-x^4
$$
because
 $$
 2p+de^2=2(d/c)+de^2\cong 2(d/c)-2(d/c)=0
 $$
 and
 $$
 p^2+de^2p-e^4=(d^2/c^2)+(d/c)d(-2/c)-(4/c^2)=-(d/c)^2-(2/c)^2\cong -1.
 $$
 This proves that $J:=(x^3y-px^5,y^2-x^4)\subseteq I.$
 We remark that
 $$
 p^2-1=(d/c)^2-1=(d^2-c^2)/c^2\cong -(2/c)^2,
 $$
 and this implies
 $$
 p^2-1\notin \n.
 $$

In order to prove that $I=J$  we need to show that the Hilbert function of $R/J$ is $1,2,2,2,1,1,1.$
We have
$$
(X^3Y,Y^2)\subseteq J^*\subseteq I^*.
$$
Further
$$
y(x^3y-px^5)-x^3(y^2-x^4)=-pyx^5+x^7\in J
$$
which implies $x^5y-(1/p)x^7 \in J.$
Thus we have
$$
x^2(x^3y-px^5)-(x^5y-(1/p) x^7 )=\frac{1-p^2}{p}x^7\in J.
$$
From this we get $x^7\in J$, hence
$$
(X^3Y,Y^2,X^7)\subseteq J^*\subseteq I^*.
$$
These ideals have the same Hilbert function so that we finally get
$$
 I=(x^3y-px^5 ,y^2-x^4)
$$
with
$$
p\notin \n,  \ \ \ \ \ \ \ p^2-1\notin \n.
$$

We have thus found three models (Case 1, Case 2a, Case 2b1) and a one dimensional family, Case 2b2.
We summarize the  models in the following table

\medskip
\begin{center}
\begin{tabular}{|l|l|l|}
  \hline
  \text{Case 1} & $I=( xy, y^4-x^6)$ &  \\  \hline
  \text{Case 2a} & $I=( x^3y, y^2-x^4)$ &  \\ \hline
  \text{Case 2b1} & $I=( x^2, x y^3-y^5)$ &  \\ \hline
  \text{Case 2b2} & $I=( x^3y- p x^5, y^2-x^4)$ & $p\notin \n$ \text{and} $p^2-1\notin \n$ \\
  \hline
\end{tabular}
\end{center}

\medskip
At this point a natural question is whether we can pass from a model to another by a
changing of generators of $\n.$

For example,   the model $I=(xy,y^4-x^6)$ of Case 1 cannot be reached by any of the other models,
because  it is quite easy to see that, however we choose the element $a\in \n,$ the ideal
$(x^3y,y^2-axy-x^4)$
does not contain the product of two minimal generators of the maximal ideal $\n.$

We are able to prove that all the models
we have found are indeed non isomorphic, but here we give a proof only for the ideals in the family of  Case 2b2.

\begin{proposition}
Let $p,q\in R$ such that  $p,q,p^2-1,q^2-1\notin \n.$ If
$\n=(x,y)=(z,v)$ and $(x^3y-px^5,y^2-x^4)=(z^3v-qz^5,v^2-z^4)$
then $p^2-q^2\in \n.$
\end{proposition}
\begin{proof}
Let $I:=(x^3y-px^5,y^2-x^4);$  we will use the equalities
$(\n/I)^3=(\overline{x}^3,\overline{x}^2\overline{y}),$
$(\n/I)^4=(\overline{x}^4),$ $(\n/I)^5=(\overline{x}^5).$

We first use  the generators $v^2-z^4$ to get $v^2\in \n^4+I\subseteq (y,x^4).$
This implies $v\in (y,x^2)$ so that $v=ex^2+by,$ with $b\notin \n.$
Since modulo $I$ we have $$v^2=e^2x^4+2ebx^2y+b^2y^2\cong e^2x^4+2ebx^2y+b^2x^4,$$
we get $2ebx^2y\in \n^4+I$ which gives $e\in \n$ and finally $$v=ax^3+by$$
with $a\in R,$ $b\notin \n.$ We also have $z=cx+dy$ with
$$det\begin{pmatrix}
    ax^2  &   c \\
   b   &  d
\end{pmatrix}=adx^2-bc\notin \n
$$
which implies $c\notin \n.$

Now, modulo $I,$ we have $ 0\cong v^2-z^4=b^2x^4-c^4x^4+t$ with
$t\in \n^5$ which implies $b^2-c^4 \in \n.$ We also have $$0\cong
z^3v-qz^5=z^3(v-qz^2)\cong c^3bpx^5-qc^5x^5+f$$ with $f\in \n^6.$
This implies $c^3bp-qc^5 \in \n, $ hence $bp-qc^2\in \n.$ Since
$b^2-c^4 \in \n$ we easily get the conclusion $p^2-q^2\in \n.$
 \end{proof}

 \bigskip
 With the methods explained before
 we can manage also the case with Hilbert function $1,3,2,1.$
 This case was the unique left case  in order  to classify,  up to isomorphism,
 Artinian Gorenstein $\res$-algebras of degree 7.
 Thus we can solve Question 4.4. of  \cite{CN07}.
 We prove that if $R/I$ is Gorenstein with Hilbert function $1,3,2,1,$ then,
 after a possible change of generators of $\n,$ either
 $$
 I=(xy,xz,yz,x^3-y^3,z^2-y^3)\ \ \ \text{ or} \ \ \ I=(x^3,y^2,yz,xz,z^2-x^2y).
 $$

\providecommand{\bysame}{\leavevmode\hbox to3em{\hrulefill}\thinspace}

\bigskip
\bigskip
\noindent
Juan Elias\\
Departament d'\`Algebra i Geometria\\
Universitat de Barcelona\\
Gran Via 585, 08007 Barcelona, Spain\\
e-mail: {\tt elias@ub.edu}

\bigskip
\noindent
Giuseppe Valla\\
Departimento di Matematica\\
Universit{\`a} di Genova\\
Via Dodecaneso 35, 16146 Genova, Italy\\
e-mail: {\tt valla@dima.unige.it}

\end{document}